\documentclass[11pt]{article}
\usepackage{amssymb}
\usepackage{graphicx}
\usepackage{amsmath}

\newtheorem{theorem}{Theorem}

\newtheorem{lemma}[theorem]{Lemma}

\newtheorem{proposition}[theorem]{Proposition}
\newtheorem{remark}[theorem]{Remark}

\newenvironment{proof}[1][Proof]{\textbf{#1.} }{\ \rule{0.5em}{0.5em}}
\input{tcilatex}

\begin{document}

\title{On elliptic boundary value problems of order $2m$ in cylindrical
domain of large size}
\author{B. Brighi - S. Guesmia \\
{\scriptsize Laboratoire Math\'{e}matiques, Informatique et Applications
(MIA)}\\
{\scriptsize 4, rue des Fr\`{e}res Lumi\`{e}re 68093 MULHOUSE CEDEX France.}%
\\
{\scriptsize E-mail: \ \ bernard.brighi@uha.fr \ \ senoussi.guesmia@uha.fr\ }%
}
\date{}
\maketitle

\begin{abstract}
We study in this work the convergence of the solution of general elliptic
boundary value problems in cylindrical domain, when some directions of the
domain go to $+\infty $.
\end{abstract}



\section{\protect\bigskip Introduction}

The present article generalizes the results of A. Rougirel and M. Chipot in 
\cite{roug}, \cite{rougirel} and \cite{ch} for the elliptic problems of
order 2. We interest in elliptic problems of order 2m, with conditions of
the Dirichlet type on cylindrical domains of $\mathbb{R}^{n}$. We study the
asymptotic behavior of the solution when the cylindrical domain becomes
unbounded in several directions. In the second section, we show under
certain conditions on data that the solution of such problems converges
towards a solution of an elliptic problem in $\mathbb{R}^{n-p}$, in Sobolev
space $H^{m}$, with a speed faster than any power of $\frac{1}{\ell }$; in
the third section, we show the same results in higher order Sobolev spaces.

Let $\omega $ be a bounded Lipschitz domain of $\mathbb{R}^{n-p}$ and $%
n>p\geqslant 1$. For a positive number $\ell $, we consider the cylinder of $%
\mathbb{R}^{n}$ 
\begin{equation*}
\Omega _{\ell }=(-\ell ,\ell )^{p}\times \omega .
\end{equation*}%
For $x=(x_{1},x_{2},\ldots ,x_{n})\in \mathbb{R}^{n}$, we will set 
\begin{equation*}
X_{1}=(x_{1},...,x_{p}),\text{ \ \ \ \ \ \ \ \ \ }X_{2}=(x_{p+1},\ldots
,x_{n}).
\end{equation*}%
We consider the boundary value problems defined by 
\begin{equation}
\left\{ 
\begin{array}{cc}
\mathcal{A}u=f & \text{in \ }\Omega _{\ell }, \\ 
\frac{\partial ^{k}u}{\partial \nu ^{k}}=0\text{ \ \ \ \ }k=0,\ldots ,m-1 & 
\text{on }\partial \Omega _{\ell },%
\end{array}%
\right.  \label{1}
\end{equation}%
\begin{equation}
\left\{ 
\begin{array}{cc}
\mathcal{A}_{\omega }u=f & \text{in \ }\omega , \\ 
\frac{\partial ^{k}u}{\partial \nu ^{k}}=0\text{ \ \ \ \ }k=0,\ldots ,m-1 & 
\text{on }\partial \omega ,%
\end{array}%
\right.  \label{2}
\end{equation}%
with

\begin{equation}
\mathcal{A}u=\sum\limits_{\left\vert \alpha \right\vert ,\left\vert \beta
\right\vert \leqslant m}(-1)^{\left\vert \alpha \right\vert }D^{\alpha
}(a_{\alpha \beta }D^{\beta }u),\text{ \ \ }\mathcal{A}_{\omega
}u=\sum\limits_{\alpha ,\beta \in N_{2}}(-1)^{\left\vert \alpha \right\vert
}D^{\alpha }(a_{\alpha \beta }D^{\beta }u),  \notag
\end{equation}%
where we have denoted by $\left\vert \alpha \right\vert $ the length of the
multi-index $\alpha ,$ $D^{\alpha }$ the partial derivative $\frac{\partial
^{\left\vert \alpha \right\vert }}{\partial ^{\alpha _{1}}x_{1}...\partial
^{\alpha _{n}}x_{n}}$, $\frac{\partial ^{k}u}{\partial \nu ^{k}}$ the $k$
derivative in the direction $\overrightarrow{\nu }$ (the unit outward normal
vector on $\partial \Omega _{\ell }$ or $\partial \omega $), and $N_{1},$ $%
N_{2}$ are given by 
\begin{equation*}
N_{1}=\left\{ \alpha \in \left\{ 0,1,\ldots,m\right\} ^{p}\times \left\{
0\right\} ^{n-p},\text{ \ }\left\vert \alpha \right\vert \leqslant m\right\}
,\ \ \ \ \ N_{2}=\left\{ \alpha \in \left\{ 0\right\} ^{p}\times \left\{
0,1,\ldots,m\right\} ^{n-p},\text{ \ }\left\vert \alpha \right\vert
\leqslant m\right\} ,
\end{equation*}%
$f$ is a function of $L^{2}(\omega )$ independent of $X_{1}$ 
\begin{equation}
f(x)=f(X_{2}),  \label{3}
\end{equation}%
the coefficients $a_{\alpha \beta }$ satisfy 
\begin{eqnarray}
a_{\alpha \beta } &\in &L^{\infty }(\mathbb{R}^{p}\times \omega )\text{ \ \
\ for \ \ }\left\vert \alpha \right\vert ,\text{ }\left\vert \beta
\right\vert \leqslant m,\text{ }  \label{4} \\
a_{\alpha \beta } &\in &C(\mathbb{R}^{p}\times \overline{\omega })\text{ \ \
\ for\ \ \ \ }\left\vert \alpha \right\vert ,\text{ }\left\vert \beta
\right\vert =m.  \label{6}
\end{eqnarray}%
and 
\begin{equation}
a_{\alpha \beta }(x)=a_{\alpha \beta }(X_{2})\text{ \ \ \ \ for \ }\alpha
\in N_{2}  \label{5}
\end{equation}%
i.e. these coefficients are only depending on the last variables $%
x_{p+1},\ldots,x_n$. We will impose the usual ellipticity condition, i.e.,
that for some $\lambda >0$ 
\begin{equation}
\sum\limits_{\left\vert \alpha \right\vert ,\left\vert \beta \right\vert
=m}(-1)^{m}a_{\alpha \beta }(x)\xi ^{\alpha +\beta }\geqslant \lambda
\left\vert \xi \right\vert ^{2m},\text{ \ a.e. \ }x\in \mathbb{R}^{p}\times
\omega \text{, \ }\forall \xi \in \mathbb{R}^{n}  \notag
\end{equation}%
where $\xi ^{\alpha +\beta }=\xi _{1}^{\alpha _{1}+\beta _{1}}\xi
_{2}^{\alpha _{2}+\beta _{2}}...\xi _{n}^{\alpha _{n}+\beta _{n}}$ ,\ \ \ $%
\left\vert \xi \right\vert ^{2}=\xi _{1}^{2}+\xi _{2}^{2}+...+\xi _{n}^{2}.$

\medskip The variational problems corresponding to (\ref{1}) and (\ref{2})
are the following 
\begin{equation}
\left\{ 
\begin{array}{c}
a(u,v):=\displaystyle\int\limits_{\Omega _{\ell }}\sum\limits_{\left\vert
\alpha \right\vert ,\left\vert \beta \right\vert \leqslant m}a_{\alpha \beta
}D^{\alpha }uD^{\beta }vdx=\int\limits_{\Omega _{\ell }}fvdx,~~\forall v\in
H_{o}^{m}(\Omega _{\ell }) \\ 
\\ 
u\in H_{o}^{m}(\Omega _{\ell }),\hfill%
\end{array}%
\right.  \label{7}
\end{equation}%
\begin{equation}
\left\{ 
\begin{array}{c}
a_{\omega }(u,v):=\displaystyle\int\limits_{\omega }\sum\limits_{\alpha
,\beta \in N_{2}}a_{\alpha \beta }D^{\alpha }uD^{\beta
}vdx=\int\limits_{\omega }fvdx,~~\forall v\in H_{o}^{m}(\omega ) \\ 
\\ 
u\in H_{o}^{m}(\omega ).\hfill%
\end{array}%
\right.  \label{8}
\end{equation}%
where $H_{o}^{m}(\Omega _{\ell })$ (resp. $H_{o}^{m}(\omega )$) is the
closure of $\mathcal{D}(\Omega _{\ell })$ (resp.$\mathcal{D}(\omega )$) in $%
H^{m}(\Omega _{\ell })$ (resp. $H^{m}(\omega )$). Then, it is well known,
see for instance \cite{lions}, that under the above assumptions, the bounded
bilinear forms $a(.,.)$ and $a_{\omega }(.,.)$ are coercive on $%
H_{o}^{m}(\Omega _{\ell })$ and $H_{o}^{m}(\omega )$ respectively, i.e.
there exist $C,C_{\ell }>0,$ $c\in \mathbb{R}$ such that 
\begin{eqnarray}
a(u,u)+c\left\Vert u\right\Vert _{L^{2}(\Omega _{\ell })}^{2} &\geqslant
&C_{\ell }\left\Vert u\right\Vert _{H^{m}(\Omega _{\ell })}^{2}\text{ \ }%
u\in H_{o}^{m}(\Omega _{\ell })\text{,}  \label{47} \\
\text{\ \ \ }a_{\omega }(u,u)+c\left\Vert u\right\Vert _{L^{2}(\omega )}^{2}
&\geqslant &C\left\Vert u\right\Vert _{H^{m}(\omega )}^{2}\text{ \ \ \ }u\in
H_{o}^{m}(\omega ),\text{\ }  \label{47'}
\end{eqnarray}%
Moreover, if we take $c=0$, there exists a unique solution $u_{\ell }$ in $%
H_{o}^{m}(\Omega _{\ell })$ to problem (\ref{7}) and a unique solution $%
u_{\infty }$ in $H_{o}^{m}(\omega )$ to problem (\ref{8}). We will also need
to assume that the constant $C_{\ell }$ in (\ref{47}) is independent of $%
\ell $, then 
\begin{equation}
a(u,u)\geqslant C\left\Vert u\right\Vert _{H^{m}(\Omega _{\ell })}^{2}\text{
\ \ \ }u\in H_{o}^{m}(\Omega _{\ell })  \label{48}
\end{equation}

\begin{remark}
We can only suppose%
\begin{equation*}
a(u,u)\geqslant \frac{C}{\ell ^{\kappa }}\left\vert u\right\vert _{m}^{2}%
\text{ \ \ \ }u\in H_{o}^{m}(\Omega _{\ell })
\end{equation*}%
where $\ 0<\kappa <1$ and $\left\vert u\right\vert
_{m}^{2}=\sum\limits_{\left\vert \alpha \right\vert \leqslant m}\displaystyle%
\int\limits_{\Omega }D^{\alpha }uD^{\alpha }vdx,$ then we have the same
results.
\end{remark}

\section{The convergence in the space $H^{m}(\Omega _{\ell _{o}})$}

We start by showing the following result :

\begin{proposition}
Let $v$ be an element of $H_{o}^{m}(\Omega _{\ell })$. Then 
\begin{equation}
v(X_{1},.)\in H_{o}^{m}(\omega )\text{ \ \ \ for almost all }X_{1}\text{ in }%
(-\ell ,\ell )^{p}.  \label{9}
\end{equation}
\end{proposition}

\begin{proof}
Using the idea of the proposition (3.1) in \cite{ch}. If we take $v\in
H_{o}^{m}(\Omega _{\ell })$. Then, there exists a sequence $\varphi _{n}$ of
element of $\mathcal{D}(\Omega _{\ell }),$ such that%
\begin{equation*}
\int\limits_{\Omega _{\ell }}\left( D^{\alpha }(v_{n}-v)\right)
^{2}dx\longrightarrow 0\text{ \ for }\left\vert \alpha \right\vert \leqslant
m.
\end{equation*}%
Thus, there exists a subsequence $v_{n^{\prime }},$ such that 
\begin{equation*}
\int\limits_{\omega }\left( D^{\alpha }(v_{n^{\prime }}-v)\left(
X_{1},.\right) \right) ^{2}dX_{2}\longrightarrow 0\text{ \ }
\end{equation*}%
for almost all $X_{1}$ in $(-\ell ,\ell )^{p}$ and for $\alpha \in N_{2}.$
Then because $v_{n^{\prime }}\left( X_{1},.\right) \in \mathcal{D}(\omega )$
for all $X_{1}$ in $(-\ell ,\ell )^{p}$ and $v_{n^{\prime }}\left(
X_{1},.\right) \longrightarrow v\left( X_{1},.\right) $ in $H_{o}^{m}(\omega
)$ for almost all $X_{1}$ in $(-\ell ,\ell )^{p},$ we have $\left( \ref{9}%
\right) .$
\end{proof}

\begin{theorem}
\label{Hm}\label{conv1}Under the assumptions (\ref{3}), (\ref{4}), (\ref{6}%
), (\ref{5}) and (\ref{48}), for all $\ell _{o}>0$ and $r>0,$ there exists a
constant $C>0$ independent of $\ell $ such that 
\begin{equation}
\left\Vert u_{\ell }-u_{\infty }\right\Vert _{H^{m}(\Omega _{\ell
_{o}})}\leqslant \frac{C}{\ell ^{r}},  \label{17}
\end{equation}%
where $u_{\ell }$ and $u_{\infty }$ are the solutions to (\ref{7}) and (\ref%
{8}) respectively.
\end{theorem}

\begin{proof}
We have 
\begin{equation*}
\int\limits_{\Omega _{\ell }}\sum\limits_{\left\vert \alpha \right\vert
,\left\vert \beta \right\vert \leqslant m}a_{\alpha \beta }D^{\alpha
}u_{\ell }D^{\beta }vdx=\int\limits_{\Omega _{\ell }}fvdx\text{ \ \ \ for
all \ }v\in H_{o}^{m}(\Omega _{\ell }),
\end{equation*}%
and also 
\begin{equation*}
\int\limits_{\omega }\sum\limits_{\alpha ,\beta \in N_{2}}a_{\alpha \beta
}D^{\alpha }u_{\infty }D^{\beta }vdx=\int\limits_{\omega }fvdx\text{\ \ \ \
\ for all \ }v\in H_{o}^{m}(\omega ).
\end{equation*}%
Applying the previous proposition, we have 
\begin{equation}
\int\limits_{\Omega _{\ell }}\sum\limits_{\left\vert \alpha \right\vert
,\left\vert \beta \right\vert \leqslant m}a_{\alpha \beta }D^{\alpha
}u_{\ell }D^{\beta }vdx=\int\limits_{(-\ell ,\ell )^{p}}\int\limits_{\omega
}\sum\limits_{\alpha ,\beta \in N_{2}}a_{\alpha \beta }D^{\alpha }u_{\infty
}D^{\beta }vdx\text{ \ \ \ \ \ \ for all \ }v\in H_{o}^{m}(\Omega _{\ell }).
\label{18}
\end{equation}%
Taking into account the independence of $u_{\infty }$ from $X_{1}$, we
obtain 
\begin{equation}
\int\limits_{\Omega _{\ell }}\sum\limits_{\left\vert \alpha \right\vert
,\left\vert \beta \right\vert \leqslant m}a_{\alpha \beta }D^{\alpha
}(u_{\ell }-u_{\infty })D^{\beta }vdx=-\int\limits_{\Omega _{\ell
}}\sum\limits_{\substack{ 0<\left\vert \alpha \right\vert ,\left\vert \beta
\right\vert \leqslant m  \\ \alpha \in N_{2},\beta \in N_{1}}}a_{\alpha
\beta }D^{\alpha }u_{\infty }D^{\beta }vdx\text{ \ \ \ for all \ \ }v\in
H_{o}^{m}(\Omega _{\ell }).  \label{10}
\end{equation}%
Using Gauss formula, and taking into account the fact that $u_{\infty }$ is
independent of $\ell $, the functions $a_{\alpha \beta }$ for $\beta \in
N_{1}$ and $\alpha \in N_{2}$ are independent of $X_{1}$, we obtain, for $%
\beta \in N_{1}$, and $\left\vert \beta \right\vert >0$ (i.e. there exists a 
$\beta _{i}\neq 0$) 
\begin{eqnarray*}
\int\limits_{\Omega _{\ell }}a_{\alpha \beta }D^{\alpha }u_{\infty }D^{\beta
}vdx &=&\int\limits_{\Omega _{\ell }}D^{\beta }\left( a_{\alpha \beta
}vD^{\alpha }u_{\infty }\right) dx \\
&=&\int\limits_{\partial \Omega _{\ell }}D^{(\beta _{1},...,\beta
_{i}-1,...,\beta _{P},0,...,0)}\left( a_{\alpha \beta }vD^{\alpha }u_{\infty
}\right) \nu _{i}dx,
\end{eqnarray*}%
then (\ref{10}) becomes 
\begin{equation}
\int\limits_{\Omega _{\ell }}\sum\limits_{\left\vert \alpha \right\vert
,\left\vert \beta \right\vert \leqslant m}a_{\alpha \beta }D^{\alpha
}(u_{\ell }-u_{\infty })D^{\beta }vdx=0\text{ \ \ \ for all \ }v\in
H_{o}^{m}(\Omega _{\ell })  \label{12}
\end{equation}%
Let $\varrho $ be a smooth function of $\mathbb{R}^{p}$, such that 
\begin{gather}
0\leqslant \varrho \leqslant 1,\text{ \ \ \ \ \ \ \ }\varrho =1\text{ on }%
\left( -\frac{1}{2},\frac{1}{2}\right) ^{p},\text{ \ \ \ \ }\varrho =0\text{
\ \ on }\mathbb{R}^{p}\backslash (-1,1)^{p},\text{\ }  \label{11} \\
\left\vert D^{\alpha }\varrho \right\vert \leqslant C\text{ \ \ \ \ }%
\left\vert \alpha \right\vert \leqslant m\text{.}  \notag
\end{gather}%
where $C$ is some constant. For $\ell _{1}\leqslant \ell $, we have 
\begin{equation*}
(u_{\ell }-u_{\infty })\varrho ^{2}\left( \frac{X_{1}}{_{\ell _{1}}}\right)
\in H_{o}^{m}(\Omega _{\ell }).
\end{equation*}%
We take in (\ref{10}) 
\begin{equation*}
v=(u_{\ell }-u_{\infty })\varrho ^{2}\left( \frac{X_{1}}{_{\ell _{1}}}%
\right) .
\end{equation*}%
We obtain 
\begin{equation}
\int\limits_{\Omega _{\ell }}\sum\limits_{\left\vert \alpha \right\vert
,\left\vert \beta \right\vert \leqslant m}a_{\alpha \beta }D^{\alpha
}(u_{\ell }-u_{\infty })D^{\beta }\left\{ (u_{\ell }-u_{\infty })\varrho
^{2}\left( \frac{X_{1}}{_{\ell _{1}}}\right) \right\} dx=0.  \label{14}
\end{equation}%
Using 
\begin{equation*}
\varrho D^{\alpha }(u_{\ell }-u_{\infty })=D^{\alpha }(\varrho (u_{\ell
}-u_{\infty }))-\sum\limits_{\alpha ^{\prime }<\alpha }\frac{1}{\ell
_{1}^{\left\vert \alpha -\alpha ^{\prime }\right\vert }}\left( _{\alpha
}^{\alpha ^{\prime }}\right) D^{\alpha ^{\prime }}(u_{\ell }-u_{\infty
})D^{\alpha -\alpha ^{\prime }}\varrho
\end{equation*}%
where\ $\left( _{\alpha }^{\alpha ^{\prime }}\right) =\left( _{\alpha
_{1}}^{\alpha _{1}^{\prime }}\right) \left( _{\alpha _{2}}^{\alpha
_{2}^{\prime }}\right) ...\left( _{\alpha _{n}}^{\alpha _{n}^{\prime
}}\right) $ and $\left( _{\alpha _{i}}^{\alpha _{i}^{\prime }}\right) =\frac{%
\alpha _{i}^{\prime }\mathbf{!}\left( \alpha _{i}-\alpha _{i}^{\prime
}\right) \mathbf{!}}{\alpha _{i}\mathbf{!}},$ the equality (\ref{14})
becomes 
\begin{multline}
\int\limits_{\Omega _{\ell _{1}}}\sum\limits_{\left\vert \alpha \right\vert
,\left\vert \beta \right\vert \leqslant m}a_{\alpha \beta }D^{\alpha
}\left\{ (u_{\ell }-u_{\infty })\varrho (\frac{X_{1}}{_{\ell _{1}}})\right\}
D^{\beta }\left\{ (u_{\ell }-u_{\infty })\varrho (\frac{X_{1}}{_{\ell _{1}}}%
)\right\} dx \\
=-\int\limits_{\Omega _{\ell _{1}}}\sum\limits_{\left\vert \alpha
\right\vert ,\left\vert \beta \right\vert \leqslant m}\sum\limits_{\beta
^{\prime }<\beta }\frac{1}{\ell _{1}^{\left\vert \beta -\beta ^{\prime
}\right\vert }}\left( _{\beta }^{\beta ^{\prime }}\right) a_{\alpha \beta
}D^{\alpha }(u_{\ell }-u_{\infty })D^{\beta ^{\prime }}((u_{\ell }-u_{\infty
})\varrho )D^{\beta -\beta ^{\prime }}\varrho dx \\
+\int\limits_{\Omega _{\ell _{1}}}\sum\limits_{\left\vert \alpha \right\vert
,\left\vert \beta \right\vert \leqslant m}\sum\limits_{\alpha ^{\prime
}<\alpha }\frac{1}{\ell _{1}^{\left\vert \alpha -\alpha ^{\prime
}\right\vert }}\left( _{\alpha }^{\alpha ^{\prime }}\right) a_{\alpha \beta
}D^{\alpha ^{\prime }}(u_{\ell }-u_{\infty })D^{\alpha -\alpha ^{\prime
}}\varrho D^{\beta }((u_{\ell }-u_{\infty })\varrho )dx.\text{ }  \label{15}
\end{multline}%
Using the Cauchy-Schwarz inequality, we can estimate all the terms of right
hand side of (\ref{15}) by 
\begin{equation}
\text{\ }\int\limits_{\Omega _{\ell _{1}}}\frac{C}{\ell _{1}^{i}}a_{\alpha
\beta }D^{\alpha }(u_{\ell }-u_{\infty })D^{\beta }((u_{\ell }-u_{\infty
})\varrho )D^{\gamma }\varrho dx\leqslant \frac{C}{\ell _{1}}\left\Vert
u_{\ell }-u_{\infty }\right\Vert _{H^{m}(\Omega _{\ell _{1}})}\left\Vert
(u_{\ell }-u_{\infty })\varrho \right\Vert _{H^{m}(\Omega _{\ell _{1}})}
\label{13}
\end{equation}%
where$\ i\geqslant 1,\ \ell _{1}\geqslant 1\ $and$\ \left\vert \alpha
\right\vert ,\left\vert \beta \right\vert ,\left\vert \gamma \right\vert
\leqslant m.$ Using the coercivity of the problem (\ref{48}) and the
estimation (\ref{13}), we obtain 
\begin{equation*}
C^{^{\prime }}\left\Vert (u_{\ell }-u_{\infty })\varrho \right\Vert
_{H^{m}(\Omega _{\ell _{1}})}^{2}\leqslant \frac{C}{\ell _{1}}\left\Vert
u_{\ell }-u_{\infty }\right\Vert _{H^{m}(\Omega _{\ell _{1}})}\left\Vert
(u_{\ell }-u_{\infty })\varrho \right\Vert _{H^{m}(\Omega _{\ell _{1}})}.
\end{equation*}%
It follows 
\begin{equation*}
\left\Vert (u_{\ell }-u_{\infty })\varrho \right\Vert _{H^{m}(\Omega _{\ell
_{1}})}\leqslant \frac{C}{\ell _{1}}\left\Vert u_{\ell }-u_{\infty
}\right\Vert _{H^{m}(\Omega _{\ell _{1}})},
\end{equation*}%
where $C$ is a constant independent of $\ell _{1}$ and $\ell .$ Since $%
\varrho =1$ on $(-\frac{1}{2},\frac{1}{2})^{p},$ we obtain 
\begin{equation*}
\left\Vert u_{\ell }-u_{\infty }\right\Vert _{H^{m}(\Omega _{\frac{\ell _{1}%
}{2}})}\leqslant \frac{C}{\ell _{1}}\left\Vert u_{\ell }-u_{\infty
}\right\Vert _{H^{m}(\Omega _{\ell _{1}})}.
\end{equation*}%
If we set $\ell _{1}=\frac{\ell }{2^{k}}$, \ $k\in \mathbb{N},$ we have 
\begin{equation*}
\left\Vert u_{\ell }-u_{\infty }\right\Vert _{H^{m}(\Omega _{\frac{\ell }{%
2^{k+1}}})}\leqslant \frac{C}{\frac{\ell }{2^{k}}}\left\Vert u_{\ell
}-u_{\infty }\right\Vert _{H^{m}(\Omega _{\frac{_{\ell }}{2^{k}}})},
\end{equation*}%
it follows that 
\begin{equation}
\left\Vert u_{\ell }-u_{\infty }\right\Vert _{H^{m}(\Omega \frac{_{\ell }}{%
2^{k}})}\leqslant \frac{C}{\ell ^{k}}\left\Vert u_{\ell }-u_{\infty
}\right\Vert _{H^{m}(\Omega _{\ell })}  \label{16}
\end{equation}%
where $C$ is only depending on $k$. Therefore, it is clear that if we can
estimate $\left\Vert u_{\ell }-u_{\infty }\right\Vert _{H^{m}(\Omega _{\ell
})}$, we have (\ref{17}).

\begin{lemma}
Under the assumption (\ref{conv1}), the following estimate holds 
\begin{equation}
\left\| u_{\ell }\right\| _{H^{m}(\Omega _{\ell })}\leqslant C\ell ^{\frac{p%
}{2}}\left\| u_{\infty }\right\| _{H^{m}(\omega )}.  \label{19}
\end{equation}
\ \ 
\end{lemma}

\begin{proof}
We set $v=u_{\ell }$ in (\ref{18}), 
\begin{equation*}
\int\limits_{\Omega _{\ell }}\sum\limits_{\left\vert \alpha \right\vert
,\left\vert \beta \right\vert \leqslant m}a_{\alpha \beta }D^{\alpha
}u_{\ell }D^{\beta }u_{\ell }dx=\int\limits_{(-\ell ,\ell
)^{p}}\int\limits_{\omega }\sum\limits_{\alpha ,\beta \in N_{2}}a_{\alpha
\beta }D^{\alpha }u_{\infty }D^{\beta }u_{\ell }dx.
\end{equation*}%
using the ellipticity of the problem (\ref{48}) and the Cauchy-Schwarz
inequality, it follows 
\begin{equation*}
C^{\prime }\left\Vert u_{\ell }\right\Vert _{H^{m}(\Omega _{\ell
})}^{2}\leqslant C\left\Vert u_{\ell }\right\Vert _{H^{m}(\Omega _{\ell
})}\left( \int\limits_{(-\ell ,\ell )^{p}}\left\Vert u_{\infty }\right\Vert
_{H^{m}(\omega )}^{2}dX_{1}\right) ^{\frac{1}{2}}.
\end{equation*}%
Then, 
\begin{equation*}
\left\Vert u_{\ell }\right\Vert _{H^{m}(\Omega _{\ell })}\leqslant C\ell ^{%
\frac{p}{2}}\left\Vert u_{\infty }\right\Vert _{H^{m}(\omega )}.
\end{equation*}%
The proof is complete.
\end{proof}

Let us come back now to the proof of the theorem. If we use (\ref{19}), the
inequality (\ref{16}) implies that 
\begin{equation*}
\left\Vert u_{\ell }-u_{\infty }\right\Vert _{H^{m}(\Omega \frac{_{\ell }}{%
2^{k}})}\leqslant \frac{C}{\ell ^{k}}\left( \left\Vert u_{\ell }\right\Vert
_{H^{m}(\Omega _{\ell })}+\left\Vert u_{\infty }\right\Vert _{H^{m}(\Omega
_{\ell })}\right) \leqslant \frac{C}{\ell ^{k}}\left( C^{{\prime }}\ell ^{%
\frac{p}{2}}\left\Vert u_{\infty }\right\Vert _{H^{m}(\omega )}+\ell ^{\frac{%
p}{2}}\left\Vert u_{\infty }\right\Vert _{H^{m}(\omega )}\right)
\end{equation*}%
from where we get 
\begin{equation*}
\left\Vert u_{\ell }-u_{\infty }\right\Vert _{H^{m}(\Omega \frac{_{\ell }}{%
2^{k}})}\leqslant \frac{C}{\ell ^{k-\frac{p}{2}}}\left\Vert u_{\infty
}\right\Vert _{H^{m}(\omega )}
\end{equation*}%
Choosing then $k$ such that $k-\frac{p}{2}>r,$ and for $\ell $ sufficient
large such that $\frac{\ell }{2^{k}}>\ell _{o},$ we obtain the desired
estimate 
\begin{equation*}
\left\Vert u_{\ell }-u_{\infty }\right\Vert _{H^{m}(\Omega _{\ell
_{o}})}\leqslant \frac{C}{\ell ^{r}},
\end{equation*}%
which completes the proof of the theorem.
\end{proof}

\section{Convergence in higher order Sobolev spaces}

In this part we suppose the functions $a_{\alpha \beta }$ verify the
following regularity conditions 
\begin{gather}
a_{\alpha \beta }\in C^{m}(\mathbb{R}^{p}\times \omega )\text{ \ \ \ for\ \
\ \ }\left| \alpha \right| ,\left| \beta \right| \leqslant m.  \label{27} \\
\left| D^{\gamma }a_{\alpha \beta }\right| \leqslant C\text{ \ \ \ on }%
\mathbb{R}^{p}\times \omega \text{ \ for \ \ }\left| \gamma \right| ,\left|
\alpha \right| ,\left| \beta \right| \leqslant m  \label{45}
\end{gather}
where $C$ is constant.

\begin{theorem}
\label{conv12}Under assumptions (\ref{3}), (\ref{5}), (\ref{48}), (\ref{27})
and (\ref{45}), then for any $\ell _{o}>0$, any $r>0,$ and any $\widetilde{%
\Omega }_{\ell _{o}}\Subset \Omega _{\ell _{o}}$ $^{(}\footnote{$^{)}$ the
closure of $\widetilde{\Omega }_{\ell _{o}}$ is in $\Omega _{\ell _{o}}.$}%
^{)}$, there exists a constant $C>0$ independent of $\ell $ such that 
\begin{eqnarray}
\left\Vert \partial _{x_{k}}\left( u_{\ell }-u_{\infty }\right) \right\Vert
_{H^{m}(\Omega _{\ell _{o}})} &\leqslant &\frac{C}{\ell ^{r}}\text{ \ \ \ \
for \ }k=1,\ldots,p  \label{22} \\
\left\Vert \partial _{x_{k}}\left( u_{\ell }-u_{\infty }\right) \right\Vert
_{H^{m}(\widetilde{\Omega }_{\ell _{o}})} &\leqslant &\frac{C}{\ell ^{r}}%
\text{ \ \ \ \ for \ }k=p+1,\ldots,n  \label{39}
\end{eqnarray}%
where $u_{\ell }$ and $u_{\infty }$ are the solutions of (\ref{7}) and (\ref%
{8}) respectively.
\end{theorem}

The idea of the proof is based on the use of finite differences, which is
possible for any type of functions, instead of derivation. Thus, for $h>0$
we define the differences of order 1 by 
\begin{equation*}
\delta _{x_{k}}v=\delta _{x_{k}}^{h}v=\frac{v(x+he_{k})-v(x)}{h}
\end{equation*}%
where $e_{k}=(0,\ldots,1,\ldots,0)$, and we define the differences of higher
order by 
\begin{eqnarray}
\delta ^{\alpha }v &=&\delta _{h}^{\alpha }v=\delta ^{\alpha _{1}}\delta
^{\alpha _{2}}\ldots\delta ^{\alpha _{n}}v  \label{24} \\
\text{where \ \ \ }\alpha &=&(\alpha _{1},\alpha _{2},\ldots,\alpha _{n})%
\text{ }\in \mathbb{N}^{n}\text{\ and \ }\delta ^{\alpha _{k}}v=\delta
_{x_{k}}\delta ^{\alpha _{k}-1}v  \notag
\end{eqnarray}%
We start by giving some properties of the finite differences analogous with
those of derivation,

\begin{lemma}
\label{par parie}Let be $\mathcal{O}$ a bounded domain in $\mathbb{R}^{n}$, $%
f\in L^{2}(\mathcal{O)}$, $\eta \in \mathcal{D(O)}$. then 
\begin{equation}
\int\limits_{\mathcal{O}}f\delta _{h}^{\alpha }\eta dx=(-1)^{\left\vert
\alpha \right\vert }\int\limits_{\mathcal{O}}\delta _{-h}^{\alpha }f\eta dx.
\label{23}
\end{equation}
\end{lemma}

\begin{proof}
Applying \cite[Lemma 3.9]{ch}, for $h$ sufficiently small, the equality (\ref%
{23}) is verified for $\left\vert \alpha \right\vert =1$. Thus, by induction
on each component of $\alpha $, and using (\ref{24}), we obtain (\ref{23}).
\end{proof}

\begin{lemma}
\label{leipnitz}Let $f$ and $g$ two functions defined on a part of $\mathbb{R%
}^{n}$. Then 
\begin{equation}
\delta ^{\alpha }(fg)=\sum\limits_{\beta \leqslant \alpha }\left( _{\alpha
}^{\beta }\right) \delta ^{\beta }f(x+(\alpha -\beta )h)\delta ^{\alpha
-\beta }g,  \label{25}
\end{equation}
\end{lemma}

\begin{proof}
This follows by induction on each component $\alpha _{k}$ of $\alpha$
\end{proof}

\begin{lemma}
Let $f$ be a function of class $C^{m}$ on the open set $\mathcal{O}\ $of $%
\mathbb{R}^{n}$. Then for any $x$ in $\mathcal{O}$,$\ h$ sufficiently small,
there exists $\xi _{h}^{x}$ of $\mathcal{O}$, such that 
\begin{equation}
\delta _{h}^{\alpha }f(x)=\frac{1}{\alpha !}D^{\alpha }f(\xi _{h}^{x})\text{
\ \ \ for\ \ \ }\left\vert \alpha \right\vert \leqslant m.  \label{46}
\end{equation}
\end{lemma}

\begin{proof}
This follows immediately from the mean-value theorem.
\end{proof}

\medskip We turn now to the proof of the theorem. Taking $w=u_{\ell
}-u_{\infty }$ in (\ref{12}), we obtain 
\begin{equation*}
\int\limits_{\Omega _{\ell }}\sum\limits_{\left\vert \alpha \right\vert
,\left\vert \beta \right\vert \leqslant m}a_{\alpha \beta }D^{\alpha
}wD^{\beta }vdx=0\text{ \ \ \ for any \ }v\in H_{o}^{m}(\Omega _{\ell }).
\end{equation*}%
For $v$ in $H_{o}^{m}(\Omega _{\ell })$ with a support disjoint of $\partial
\left( -\ell ,\ell \right) ^{p}\times \overline{\omega }$ if $\gamma \in
N_{1}$ and with a support in $\Omega _{\ell }$ if $\gamma \notin N_{1}$, and 
$h$ sufficiently small, we can replace in the above expression $v$ by $%
(-1)^{\left\vert \gamma \right\vert }\delta _{-h}^{\gamma }v$, with $%
\left\vert \gamma \right\vert \leqslant m$. Using the permutation of the
derivation and the finite difference, and the lemma \ref{par parie}, we
obtain 
\begin{equation*}
\int\limits_{\Omega _{\ell }}\sum\limits_{\left\vert \alpha \right\vert
,\left\vert \beta \right\vert \leqslant m}\delta _{h}^{\gamma }(a_{\alpha
\beta }D^{\alpha }w)D^{\beta }vdx=0.
\end{equation*}%
The lemma \ref{leipnitz} with $f=a_{\alpha \beta }$ \ and $g=D^{\alpha }w$
gives 
\begin{multline}
\int\limits_{\Omega _{\ell }}\sum\limits_{\left\vert \alpha \right\vert
,\left\vert \beta \right\vert \leqslant m}a_{\alpha \beta }(x+\gamma
h)\delta ^{\gamma }D^{\alpha }wD^{\beta }vdx \\
=-\int\limits_{\Omega _{\ell }}\sum\limits_{\left\vert \alpha \right\vert
,\left\vert \beta \right\vert \leqslant m}\sum\limits_{0<\sigma \leqslant
\gamma }\left( _{\gamma }^{\sigma }\right) \delta ^{\sigma }a_{\alpha \beta
}(x+\left( \gamma -\sigma \right) h)\delta ^{\gamma -\sigma }D^{\alpha
}wD^{\beta }vdx.  \label{26}
\end{multline}%
Let $\ell _{o}<\ell _{1}\leqslant \ell ,$ and $\Omega _{\ell _{o}}^{\prime
},\Omega _{\ell _{1}}^{\prime }$ two bounded domain of $\mathbb{R}^{n}$,
such that $\Omega _{\ell _{o}}^{\prime }=\Omega _{\ell _{o}}$, $\Omega
_{\ell _{1}}^{\prime }=\Omega _{\ell ^{\prime }}$ with $\ell _{o}<\ell
^{\prime }<\ell _{1}$ if $\gamma \in N_{1}$, and $\Omega _{\ell
_{o}}^{\prime }=\widetilde{\Omega }_{\ell _{o}}\Subset \Omega _{\ell
_{1}}^{\prime }\Subset \Omega _{\ell _{1}}$ if $\gamma \notin N_{1}$. Let us
denote by $\varrho $ a smooth function with compact support in $\left( -\ell
^{\prime },\ell ^{\prime }\right) ^{p}\times \overline{\omega }$ if $\gamma
\in N_{1},$ and with compact support in $\Omega _{\ell _{1}}^{\prime }$ if $%
\gamma \notin N_{1}$, such that in both cases we have 
\begin{equation*}
0\leqslant \varrho \leqslant 1,\text{ \ \ \ \ }\varrho =1\text{ \ \ \ on \ }%
\Omega _{\ell _{o}}^{\prime }.
\end{equation*}%
We take $v=\delta ^{\gamma }w\varrho ^{4m}$ in (\ref{26}) for $h$ small
enough , and using equalities 
\begin{eqnarray}
D^{\beta }\left\{ (\delta ^{\gamma }w\varrho ^{2m})\varrho ^{2m}\right\}
&=&\varrho ^{2m}D^{\beta }(\delta ^{\gamma }w\varrho
^{2m})+\sum\limits_{\beta ^{{\prime }}<\beta }\left( _{\beta }^{\beta ^{{%
\prime }}}\right) D^{\beta ^{{\prime }}}(\delta ^{\gamma }w\varrho
^{2m})D^{\beta -\beta ^{{\prime }}}\varrho ^{2m}  \notag \\
\varrho ^{2m}D^{\alpha }\delta ^{\gamma }w &=&D^{\alpha }(\varrho
^{2m}\delta ^{\gamma }w)-\sum\limits_{\alpha ^{{\prime }}<\alpha }\left(
_{\alpha }^{\alpha ^{\prime }}\right) D^{\alpha ^{{\prime }}}\delta ^{\gamma
}wD^{\alpha -\alpha ^{{\prime }}}\varrho ^{2m},  \label{28}
\end{eqnarray}%
we see that (\ref{26}) becomes 
\begin{gather}
\int\limits_{\Omega _{\ell _{1}}^{\prime }}\sum\limits_{\left\vert \alpha
\right\vert ,\left\vert \beta \right\vert \leqslant m}a_{\alpha \beta
}(x+\gamma h)D^{\alpha }\left( \delta ^{\gamma }w\varrho ^{2m}\right)
D^{\beta }\left( \delta ^{\gamma }w\varrho ^{2m}\right) dx=  \notag \\
-\int\limits_{\Omega _{\ell _{1}}^{\prime }}\sum\limits_{\left\vert \alpha
\right\vert ,\left\vert \beta \right\vert \leqslant m}\sum\limits_{0<\sigma
\leqslant \gamma }\sum\limits_{\beta ^{{\prime }}\leqslant \beta }\left(
_{\beta }^{\beta ^{{\prime }}}\right) \left( _{\gamma }^{\sigma }\right)
\delta ^{\sigma }a_{\alpha \beta }(x+\left( \gamma -\sigma \right)
h)D^{\alpha }\delta ^{\gamma -\sigma }wD^{\beta ^{{\prime }}}(\delta
^{\gamma }w\varrho ^{2m})D^{\beta -\beta ^{{\prime }}}\varrho ^{2m}dx  \notag
\\
-\int\limits_{\Omega _{\ell _{1}}^{\prime }}\sum\limits_{\left\vert \alpha
\right\vert ,\left\vert \beta \right\vert \leqslant m}\sum\limits_{\beta ^{{%
\prime }}<\beta }\left( _{\beta }^{\beta ^{\prime }}\right) a_{\alpha \beta
}(x+\gamma h)D^{\alpha }\delta ^{\gamma }wD^{\beta ^{{\prime }}}(\delta
^{\gamma }w\varrho ^{2m})D^{\beta -\beta ^{{\prime }}}\varrho ^{2m}dx  \notag
\\
+\int\limits_{\Omega _{\ell _{1}}^{\prime }}\sum\limits_{\left\vert \alpha
\right\vert ,\left\vert \beta \right\vert \leqslant m}\sum\limits_{\alpha ^{{%
\prime }}<\alpha }\left( _{\alpha }^{\alpha ^{\prime }}\right) a_{\alpha
\beta }(x+\gamma h)D^{\alpha ^{{\prime }}}\delta ^{\gamma }wD^{\beta }\left(
\delta ^{\gamma }w\varrho ^{2m}\right) D^{\alpha -\alpha ^{{\prime }%
}}\varrho ^{2m}dx.  \label{3A}
\end{gather}%
We estimate one by one the three terms of the right hand side. The first
term is the sum of terms of the form 
\begin{equation*}
\int\limits_{\Omega _{\ell _{1}}^{\prime }}C\delta ^{\sigma }a_{\alpha \beta
}(x+\left( \gamma -\sigma \right) h)D^{\alpha }\delta ^{\gamma -\sigma
}wD^{\beta ^{{\prime }}}(\delta ^{\gamma }w\varrho ^{2m})D^{\beta -\beta ^{{%
\prime }}}\varrho ^{2m}dx
\end{equation*}%
such that $C$ is a constant, $0<\sigma \leqslant \gamma $, $\beta ^{{\prime }%
}\leqslant \beta $ and $\left\vert \alpha \right\vert ,\left\vert \beta
\right\vert \leqslant m$. Using (\ref{46}) and the fact that the function $%
\varrho $ and these derivatives are bounded, and the Cauchy-Schwarz
inequality, we can estimate these terms 
\begin{multline}
\left\vert ~\int\limits_{\Omega _{\ell _{1}}^{\prime
}}\sum\limits_{\left\vert \alpha \right\vert ,\left\vert \beta \right\vert
\leqslant m}\sum\limits_{0<\sigma \leqslant \gamma }\sum\limits_{\beta ^{{%
\prime }}\leqslant \beta }\left( _{\beta }^{\beta ^{{\prime }}}\right)
\left( _{\gamma }^{\sigma }\right) \delta ^{\sigma }a_{\alpha \beta
}(x+\left( \gamma -\sigma \right) h)D^{\alpha }\delta ^{\gamma -\sigma
}wD^{\beta ^{{\prime }}}(\delta ^{\gamma }w\varrho ^{2m})D^{\beta -\beta ^{{%
\prime }}}\varrho ^{2m}dx~\right\vert \\
\leqslant C\sum\limits_{\sigma <\gamma }\left\Vert \delta ^{\gamma }w\varrho
^{2m}\right\Vert _{H^{m}(\Omega _{\ell _{1}}^{\prime })}\left\Vert \delta
^{\sigma }w\right\Vert _{H^{m}(\Omega _{\ell _{1}}^{\prime })}.  \label{34}
\end{multline}%
The third term is the sum of terms of the form 
\begin{equation*}
\int\limits_{\Omega _{\ell _{1}}^{\prime }}Ca_{\alpha \beta }(x+\gamma
h)D^{\alpha ^{{\prime }}}\delta ^{\gamma }wD^{\beta }\left( \delta ^{\gamma
}w\varrho ^{2m}\right) D^{\alpha -\alpha ^{{\prime }}}\varrho ^{2m}dx
\end{equation*}%
where $C$ is a constant, $\alpha ^{{\prime }}<$ $\alpha $ and $\left\vert
\alpha \right\vert ,\left\vert \beta \right\vert \leqslant m$. Using (\ref%
{27}), (\ref{45}) and the fact that the function $\varrho $ and these
derivatives are bounded, and the Cauchy-Schwarz inequality, we obtain 
\begin{multline}
\left\vert ~\int\limits_{\Omega _{\ell _{1}}^{\prime
}}\sum\limits_{\left\vert \alpha \right\vert ,\left\vert \beta \right\vert
\leqslant m}\sum\limits_{\beta ^{\prime }<\beta }\left( _{\beta }^{\beta
^{\prime }}\right) a_{\alpha \beta }(x+\gamma h)D^{\alpha }\delta ^{\gamma
}wD^{\beta ^{\prime }}(\delta ^{\gamma }w\varrho ^{2m})D^{\beta -\beta
^{\prime }}\varrho ^{2m}dx~\right\vert \\
\leqslant C\left\Vert \delta ^{\gamma }w\right\Vert _{H^{m-1}(\Omega _{\ell
_{1}}^{\prime })}\left\Vert \delta ^{\gamma }w\varrho ^{2m}\right\Vert
_{H^{m}(\Omega _{\ell _{1}}^{\prime })}.  \label{35}
\end{multline}%
For the second term, a direct estimate as for the other terms is not
sufficient. We first write%
\begin{equation*}
D^{\beta ^{\prime }}(\delta ^{\gamma }w\varrho ^{2m})=\sum\limits_{\beta
^{\prime \prime }\leqslant \beta ^{\prime }}\left( _{\beta ^{\prime
}}^{\beta ^{\prime \prime }}\right) D^{\beta ^{\prime \prime }}\delta
^{\gamma }wD^{\beta ^{\prime }-\beta ^{\prime \prime }}\varrho ^{2m},
\end{equation*}%
and for $\left\vert \tau \right\vert \leqslant m$%
\begin{equation*}
D^{\tau }\varrho ^{2m}=\varrho ^{m}\psi _{\varrho },
\end{equation*}%
where $\psi _{\varrho }$ is a sum and product of $\varrho $ and these
derivatives. Then 
\begin{multline*}
\int\limits_{\Omega _{\ell _{1}}^{\prime }}\sum\limits_{\left\vert \alpha
\right\vert ,\left\vert \beta \right\vert \leqslant m}\sum\limits_{\beta ^{{%
\prime }}<\beta }\left( _{\beta }^{\beta ^{{\prime }}}\right) a_{\alpha
\beta }(x+\gamma h)D^{\alpha }\delta ^{\gamma }wD^{\beta ^{{\prime }%
}}(\delta ^{\gamma }w\varrho ^{2m})D^{\beta -\beta ^{{\prime }}}\varrho
^{2m}dx \\
=\int\limits_{\Omega _{\ell _{1}}^{\prime }}\sum\limits_{\left\vert \alpha
\right\vert ,\left\vert \beta \right\vert \leqslant m}\sum\limits_{\beta ^{{%
\prime }}<\beta }\sum\limits_{\beta ^{\prime \prime }\leqslant \beta
^{\prime }}\left( _{\beta ^{\prime }}^{\beta ^{\prime \prime }}\right)
\left( _{\beta }^{\beta ^{\prime }}\right) a_{\alpha \beta }(x+\gamma
h)D^{\alpha }\delta ^{\gamma }w\varrho ^{2m}D^{\beta ^{\prime \prime
}}\delta ^{\gamma }w\psi _{\varrho }dx.
\end{multline*}%
Using (\ref{28}), we obtain 
\begin{multline}
\int\limits_{\Omega _{\ell _{1}}^{\prime }}\sum\limits_{\left\vert \alpha
\right\vert ,\left\vert \beta \right\vert \leqslant m}\sum\limits_{\beta ^{{%
\prime }}<\beta }\left( _{\beta }^{\beta ^{{\prime }}}\right) a_{\alpha
\beta }(x+\gamma h)D^{\alpha }\delta ^{\gamma }wD^{\beta ^{{\prime }%
}}(\delta ^{\gamma }w\varrho ^{2m})D^{\beta -\beta ^{{\prime }}}\varrho
^{2m}dx \\
=\int\limits_{\Omega _{\ell _{1}}^{\prime }}\sum\limits_{\left\vert \alpha
\right\vert ,\left\vert \beta \right\vert \leqslant m}\sum\limits_{\beta ^{{%
\prime }}<\beta }\sum\limits_{\beta ^{\prime \prime }\leqslant \beta
^{\prime }}\left( _{\beta ^{\prime }}^{\beta ^{\prime \prime }}\right)
\left( _{\beta }^{\beta ^{{\prime }}}\right) a_{\alpha \beta }(x+\gamma
h)D^{\alpha }(\varrho ^{2m}\delta ^{\gamma }w)D^{\beta ^{\prime \prime
}}\delta ^{\gamma }w\psi _{\varrho }dx \\
-\int\limits_{\Omega _{\ell _{1}}^{\prime }}\sum\limits_{\left\vert \alpha
\right\vert ,\left\vert \beta \right\vert \leqslant m}\sum\limits_{\beta ^{{%
\prime }}<\beta }\sum\limits_{\beta ^{\prime \prime }\leqslant \beta
^{\prime }}\sum\limits_{\alpha ^{{\prime }}<\alpha }\left( _{\alpha
}^{\alpha ^{\prime }}\right) \left( _{\beta }^{\beta ^{\prime {\prime }%
}}\right) \left( _{\beta }^{\beta ^{{\prime }}}\right) a_{\alpha \beta
}(x+\gamma h)D^{\alpha ^{{\prime }}}\delta ^{\gamma }wD^{\beta ^{\prime
\prime }}\delta ^{\gamma }w\psi _{\varrho }D^{\alpha -\alpha ^{{\prime }%
}}\varrho ^{2m}dx.  \label{29}
\end{multline}%
Therefore, the second term is a sum of two terms. The first term is a sum of
terms of the form 
\begin{equation*}
\int\limits_{\Omega _{\ell _{1}}^{\prime }}Ca_{\alpha \beta }(x+\gamma
h)D^{\alpha }(\varrho ^{2m}\delta ^{\gamma }w)D^{\beta }\delta ^{\gamma
}w\psi _{\varrho }dx
\end{equation*}%
where $\left\vert \alpha \right\vert \leqslant m$, $\left\vert \beta
\right\vert <m$. Using (\ref{27}), (\ref{45}) and the fact that the function 
$\varrho $ and these derivatives are bounded, and the Cauchy-Schwarz
inequality, we obtain 
\begin{equation}
\left\vert ~\int\limits_{\Omega _{\ell _{1}}^{\prime }}Ca_{\alpha \beta
}(x+\gamma h)D^{\alpha }(\varrho ^{2m}\delta ^{\gamma }w)D^{\beta }\delta
^{\gamma }w\psi _{\varrho }dx~\right\vert \leqslant C\left\Vert \delta
^{\gamma }w\right\Vert _{H^{m-1}(\Omega _{\ell _{1}}^{\prime })}\left\Vert
\delta ^{\gamma }w\varrho ^{2m}\right\Vert _{H^{m}(\Omega _{\ell
_{1}}^{\prime })}.  \label{30}
\end{equation}%
The second term can be developed as a sum of terms of the form 
\begin{equation*}
\int\limits_{\Omega _{\ell _{1}}^{\prime }}Ca_{\alpha \beta }(x+\gamma
h)D^{\alpha }\delta ^{\gamma }wD^{\beta }\delta ^{\gamma }w\psi _{\varrho
}D^{\alpha ^{\prime }}\varrho ^{2m}dx
\end{equation*}%
where $\left\vert \alpha ^{\prime }\right\vert \leqslant m$, $\left\vert
\alpha \right\vert ,\left\vert \beta \right\vert <m$, and again using (\ref%
{27}), (\ref{45}) and the fact that the function $\varrho $ and these
derivatives are bounded, and the Cauchy-Schwarz inequality, we obtain 
\begin{equation}
\left\vert ~\int\limits_{\Omega _{\ell _{1}}^{\prime }}Ca_{\alpha \beta
}(x+\gamma h)D^{\alpha ^{{\prime }}}\delta ^{\gamma }wD^{\beta ^{\prime
\prime }}\delta ^{\gamma }w\psi _{\varrho }D^{\alpha -\alpha ^{{\prime }%
}}\varrho ^{2m}dx~\right\vert \leqslant C\left\Vert \delta ^{\gamma
}w\right\Vert _{H^{m-1}(\Omega _{\ell _{1}}^{\prime })}^{2}.  \label{31}
\end{equation}%
By (\ref{30}) and (\ref{31}) we can estimate the second term of the second
member of (\ref{3A}) 
\begin{multline}
\left\vert ~\int\limits_{\Omega _{\ell _{1}}^{\prime
}}\sum\limits_{\left\vert \alpha \right\vert ,\left\vert \beta \right\vert
\leqslant m}\sum\limits_{\beta ^{{\prime }}<\beta }\left( _{\beta }^{\beta ^{%
{\prime }}}\right) a_{\alpha \beta }(x+\gamma h)D^{\alpha }\delta ^{\gamma
}wD^{\beta ^{{\prime }}}(\delta ^{\gamma }w\varrho ^{2m})D^{\beta -\beta ^{{%
\prime }}}\varrho ^{2m}dx~\right\vert \\
\leqslant C\left\Vert \delta ^{\gamma }w\right\Vert _{H^{m-1}(\Omega _{\ell
_{1}}^{\prime })}^{2}+C^{\prime }\left\Vert \delta ^{\gamma }w\right\Vert
_{H^{m-1}(\Omega _{\ell _{1}})}\left\Vert \delta ^{\gamma }w\varrho
^{2m}\right\Vert _{H^{m}(\Omega _{\ell _{1}}^{\prime })},  \label{33}
\end{multline}%
and finally by\ (\ref{34}), (\ref{35}) and (\ref{33}), we find the desired
estimate%
\begin{multline*}
\int\limits_{\Omega _{\ell _{1}}^{\prime }}\sum\limits_{\left\vert \alpha
\right\vert ,\left\vert \beta \right\vert \leqslant m}a_{\alpha \beta
}(x+\gamma h)\delta ^{\gamma }D^{\alpha }\left( \delta ^{\gamma }w\varrho
^{2m}\right) D^{\beta }\left( \delta ^{\gamma }w\varrho ^{2m}\right) dx \\
\leqslant C_{1}\sum\limits_{\sigma <\gamma }\left\Vert \delta ^{\gamma
}w\varrho ^{2m}\right\Vert _{H^{m}(\Omega _{\ell _{1}}^{\prime })}\left\Vert
\delta ^{\sigma }w\right\Vert _{H^{m}(\Omega _{\ell _{1}}^{\prime
})}+C_{2}\left\Vert \delta ^{\gamma }w\right\Vert _{H^{m-1}(\Omega _{\ell
_{1}}^{\prime })}\left\Vert \delta ^{\gamma }w\varrho ^{2m}\right\Vert
_{H^{m}(\Omega _{\ell _{1}}^{\prime })} \\
+C_{3}\left\Vert \delta ^{\gamma }w\right\Vert _{H^{m-1}(\Omega _{\ell
_{1}}^{\prime })}^{2}+C_{4}\left\Vert \delta ^{\gamma }w\right\Vert
_{H^{m-1}(\Omega _{\ell _{1}}^{\prime })}\left\Vert \delta ^{\gamma
}w\varrho ^{2m}\right\Vert _{H^{m}(\Omega _{\ell _{1}}^{\prime })}.
\end{multline*}%
Using the coercivity of the problem (\ref{48}) and the Young inequality, 
it follows that 
\begin{equation*}
\left\Vert \delta ^{\gamma }w\varrho ^{2m}\right\Vert _{H^{m}(\Omega _{\ell
_{1}}^{\prime })}^{2}\leqslant C_{\varepsilon }\left( \sum\limits_{\sigma
<\gamma }\left\Vert \delta ^{\sigma }w\right\Vert _{H^{m}(\Omega _{\ell
_{1}}^{\prime })}^{2}+\left\Vert \delta ^{\gamma }w\right\Vert
_{H^{m-1}(\Omega _{\ell _{1}}^{\prime })}^{2}\right) +\varepsilon
C\left\Vert \delta ^{\gamma }w\varrho ^{2m}\right\Vert _{H^{m}(\Omega _{\ell
_{1}}^{\prime })}^{2}.
\end{equation*}%
Taking $\varepsilon =\frac{1}{2C}$, we obtain 
\begin{equation*}
\left\Vert \delta ^{\gamma }w\varrho ^{2m}\right\Vert _{H^{m}(\Omega _{\ell
_{1}}^{\prime })}^{2}\leqslant C\left( \sum\limits_{\sigma <\gamma
}\left\Vert \delta ^{\sigma }w\right\Vert _{H^{m}(\Omega _{\ell
_{1}}^{\prime })}^{2}+\left\Vert \delta ^{\gamma }w\right\Vert
_{H^{m-1}(\Omega _{\ell _{1}}^{\prime })}^{2}\right) .
\end{equation*}%
Since $\varrho =1$ on $\ \Omega _{\ell _{o}}^{\prime },$ we have 
\begin{equation}
\left\Vert \delta ^{\gamma }w\right\Vert _{H^{m}(\Omega _{\ell _{o}}^{\prime
})}^{2}\leqslant C\left( \sum\limits_{\sigma <\gamma }\left\Vert \delta
^{\sigma }w\right\Vert _{H^{m}(\Omega _{\ell _{1}}^{\prime
})}^{2}+\left\Vert \delta ^{\gamma }w\right\Vert _{H^{m-1}(\Omega _{\ell
_{1}}^{\prime })}^{2}\right) .  \label{40}
\end{equation}%
If $\left\vert \gamma \right\vert =1$ (i.e. $\gamma =e_{k}$) then for $%
\sigma =0$ in (\ref{40}) we get 
\begin{equation*}
\left\Vert \delta _{x_{k}}w\right\Vert _{H^{m}(\Omega _{\ell _{o}}^{\prime
})}^{2}\leqslant C\left( \left\Vert w\right\Vert _{H^{m}(\Omega _{\ell
_{1}}^{\prime })}^{2}+\left\Vert \delta _{x_{k}}w\right\Vert
_{H^{m-1}(\Omega _{\ell _{1}}^{\prime })}^{2}\right) \text{ \ for }k=1,...,n.
\end{equation*}%
For another bounded domain verifying the same conditions as $\Omega _{\ell
_{1}}^{\prime }$ and containing the closure of $\Omega _{\ell _{1}}^{\prime
} $ (we still denote it by $\Omega _{\ell _{1}}^{\prime }$), and using \cite[%
Lemma 3.10]{ch}, we obtain 
\begin{equation*}
\left\Vert \delta _{x_{k}}w\right\Vert _{H^{m}(\Omega _{\ell _{o}}^{\prime
})}^{2}\leqslant C\left\Vert w\right\Vert _{H^{m}(\Omega _{\ell
_{1}}^{\prime })}^{2}\leqslant C\left\Vert w\right\Vert _{H^{m}(\Omega
_{\ell _{1}})}^{2}\text{ \ for }k=1,...,n.
\end{equation*}%
For fixed $\ell _{1}$, and by theorem \ref{Hm}, it holds 
\begin{equation*}
\left\Vert \delta _{x_{k}}^{h}w\right\Vert _{H^{m}(\Omega _{\ell
_{o}}^{\prime })}^{2}\leqslant \frac{C}{\ell ^{2r}}\text{ \ for }k=1,...,n,
\end{equation*}%
and thus%
\begin{equation*}
\left\Vert \partial ^{\alpha }\delta _{x_{k}}^{h}w\right\Vert _{L^{2}(\Omega
_{\ell _{o}}^{\prime })}^{2}\leqslant \frac{C}{\ell ^{2r}}\text{ \ \ \ for }%
\left\vert \alpha \right\vert =m,\text{ \ }k=1,...,n.
\end{equation*}%
Then the sequence $\left( \partial ^{\alpha }\delta _{x_{k}}^{h}w\right)
_{h} $ is bounded in $L^{2}(\Omega _{\ell _{o}}^{\prime })$ and we can
extract a subsequence $\left( \partial ^{\alpha }\delta
_{x_{k}}^{h_{n}}w\right) _{n\in \mathbb{N}}$ $\left( h_{n}\longrightarrow
0\right) $ which converges weakly in $L^{2}(\Omega _{\ell _{o}}^{\prime })$
to some function $w_{\alpha ,k}$ of $L^{2}(\Omega _{\ell _{o}}^{\prime }).$
We then obtain 
\begin{align*}
\left\Vert w_{\alpha ,k}\right\Vert _{L^{2}(\Omega _{\ell _{o}}^{\prime
})}^{2}& =\left\langle w_{\alpha ,k},w_{\alpha ,k}\right\rangle
_{L^{2}(\Omega _{\ell _{o}}^{\prime })}=\lim\limits_{n\longrightarrow
0}\left\langle w_{\alpha ,k},\partial ^{\alpha }\delta
_{x_{k}}^{h_{n}}w\right\rangle _{L^{2}(\Omega _{\ell _{o}}^{\prime })} \\
& \leqslant \lim\limits_{n\longrightarrow 0}\left\Vert \partial ^{\alpha
}\delta _{x_{k}}^{h_{n}}w\right\Vert _{L^{2}(\Omega _{\ell _{o}}^{\prime
})}\left\Vert w_{\alpha ,k}\right\Vert _{L^{2}(\Omega _{\ell _{o}}^{\prime
})}\leqslant \frac{C}{\ell ^{2r}}\left\Vert w_{\alpha ,k}\right\Vert
_{L^{2}(\Omega _{\ell _{o}}^{\prime })}.
\end{align*}%
It follows 
\begin{equation}
\left\Vert w_{\alpha ,k}\right\Vert _{L^{2}(\Omega _{\ell _{o}}^{\prime
})}\leqslant \frac{C}{\ell ^{2r}}.  \label{41}
\end{equation}%
In other way 
\begin{eqnarray*}
\partial ^{\alpha }\delta _{x_{k}}^{h_{n}}w &\longrightarrow &\partial
^{\alpha }\partial _{x_{k}}w\text{ \ \ in }\mathcal{D}^{\prime }(\Omega
_{\ell _{o}}^{\prime })\text{ \ } \\
\text{ \ \ \ }\partial ^{\alpha }\delta _{x_{k}}^{h_{n}}w &\longrightarrow
&w_{\alpha ,k}\text{ \ \ \ \ \ \ in }\mathcal{D}^{\prime }(\Omega _{\ell
_{o}}^{\prime }),
\end{eqnarray*}%
and by uniqueness of the limit, we deduce that $\partial ^{\alpha }\partial
_{x_{k}}w=w_{\alpha ,k}\in L^{2}(\Omega _{\ell _{o}}^{\prime })$, and the
proof is completed by (\ref{41}).

\begin{theorem}
Under the assumptions (\ref{3}), (\ref{5}), (\ref{48}), (\ref{27}) and (\ref%
{45}), then for any $\ell _{o}>0$, any $r>0$ and $\Omega _{\ell
_{o}}^{\prime }\Subset \Omega _{\ell _{o}}$, we have $u_{\ell }-u_{\infty
}\in H^{2m}(\Omega _{\ell _{o}}^{\prime })$, and there exists a constant $%
C>0 $ independent of $\ell $ such that 
\begin{equation}
\left\Vert u_{\ell }-u_{\infty }\right\Vert _{H^{2m}(\Omega _{\ell
_{o}}^{\prime })}\leqslant \frac{C}{\ell ^{r}}.  \label{42}
\end{equation}
\end{theorem}

\begin{proof}
It is enough to show 
\begin{equation}
\left\Vert D^{\alpha }\left( u_{\ell }-u_{\infty }\right) \right\Vert
_{H^{m}(\Omega _{\ell _{o}}^{\prime })}\leqslant \frac{C}{\ell ^{r}}\text{ \
for }\left\vert \alpha \right\vert \leqslant m,  \label{43}
\end{equation}%
by the theorem \ref{conv12}, the inequality (\ref{43}) is verified for $%
\left\vert \alpha \right\vert =1.$ We show the result by induction. Let us
suppose that for $\left\vert \sigma \right\vert <\left\vert \alpha
\right\vert \leqslant m$ we have 
\begin{equation}
\left\Vert D^{\sigma }\left( u_{\ell }-u_{\infty }\right) \right\Vert
_{H^{m}(\Omega _{\ell _{o}}^{\prime \prime })}\leqslant \frac{C}{\ell ^{r}}%
\text{,}  \label{44}
\end{equation}%
for any open $\Omega _{\ell _{o}}^{\prime \prime }$such that $\Omega _{\ell
_{o}}^{\prime }\Subset \Phi \Subset \Omega _{\ell _{o}}^{\prime \prime
}\Subset \Omega _{\ell _{o}}$. Using (\ref{40}), we obtain for $h$ small
enough%
\begin{equation*}
\left\Vert \delta ^{\alpha }\left( u_{\ell }-u_{\infty }\right) \right\Vert
_{H^{m}(\Omega _{\ell _{o}}^{\prime })}^{2}\leqslant C\left(
\sum\limits_{\sigma <\alpha }\left\Vert \delta ^{\sigma }\left( u_{\ell
}-u_{\infty }\right) \right\Vert _{H^{m}(\Phi )}^{2}+\left\Vert \delta
^{\alpha }\left( u_{\ell }-u_{\infty }\right) \right\Vert _{H^{m-1}(\Phi
)}^{2}\right) ,
\end{equation*}%
and using also \cite[Lemma 3.10]{ch} several times, we get 
\begin{eqnarray*}
\left\Vert \delta ^{\alpha }\left( u_{\ell }-u_{\infty }\right) \right\Vert
_{H^{m}(\Omega _{\ell _{o}}^{\prime })}^{2} &\leqslant &C\left(
\sum\limits_{\sigma <\alpha }\left\Vert D^{\sigma }\left( u_{\ell
}-u_{\infty }\right) \right\Vert _{H^{m}(\Omega _{\ell _{o}}^{\prime \prime
})}^{2}+\left\Vert D^{\alpha }\left( u_{\ell }-u_{\infty }\right)
\right\Vert _{H^{m-1}(\Omega _{\ell _{o}}^{\prime \prime })}^{2}\right)  \\
&\leqslant &C\sum\limits_{\sigma <\alpha }\left\Vert D^{\sigma }\left(
u_{\ell }-u_{\infty }\right) \right\Vert _{H^{m}(\Omega _{\ell _{o}}^{\prime
\prime })}^{2}.
\end{eqnarray*}%
Thanks to (\ref{44}), we obtain 
\begin{equation*}
\left\Vert \delta ^{\alpha }\left( u_{\ell }-u_{\infty }\right) \right\Vert
_{H^{m}(\Omega _{\ell _{o}}^{\prime })}^{2}\leqslant \frac{C}{\ell ^{r}},
\end{equation*}%
it holds that 
\begin{equation*}
\left\Vert D^{\beta }\delta ^{\alpha }\left( u_{\ell }-u_{\infty }\right)
\right\Vert _{L^{2}(\Omega _{\ell _{o}}^{\prime })}^{2}\leqslant \frac{C}{%
\ell ^{2r}}\text{ \ \ for }\left\vert \beta \right\vert \leqslant m.
\end{equation*}%
The sequence $\left( D^{\beta }\delta ^{\alpha }\left( u_{\ell }-u_{\infty
}\right) \right) _{h>0}$ is bounded in $L^{2}(\Omega _{\ell _{o}}^{\prime })$%
, and we can find a subsequence \newline
$\left( D^{\beta }\delta ^{\alpha }\left( u_{\ell }-u_{\infty }\right)
\right) _{h_{n}}$ $\left( h_{n}\longrightarrow 0\right) $ converging weakly
in $L^{2}(\Omega _{\ell _{o}}^{\prime })$ to a function $w_{\alpha ,\beta ,k}
$ of $L^{2}(\Omega _{\ell _{o}}^{\prime }).$ Then, we have 
\begin{gather*}
\left\Vert w_{\alpha ,\beta ,k}\right\Vert _{L^{2}(\Omega _{\ell
_{o}}^{\prime })}^{2}=\left\langle w_{\alpha ,\beta ,k},w_{\alpha ,\beta
,k}\right\rangle _{L^{2}(\Omega _{\ell _{o}}^{\prime })}=\lim \left\langle
w_{\alpha ,k},D^{\beta }\delta ^{\alpha }\left( u_{\ell }-u_{\infty }\right)
\right\rangle _{L^{2}(\Omega _{\ell _{o}}^{\prime })} \\
\leqslant \liminf \left\Vert D^{\beta }\delta ^{\alpha }\left( u_{\ell
}-u_{\infty }\right) \right\Vert _{L^{2}(\Omega _{\ell _{o}}^{\prime
})}\left\Vert w_{\alpha ,\beta ,k}\right\Vert _{L^{2}(\Omega _{\ell
_{o}}^{\prime })}\leqslant \frac{C}{\ell ^{2r}}\left\Vert w_{\alpha ,\beta
,k}\right\Vert _{L^{2}(\Omega _{\ell _{o}}^{\prime })}
\end{gather*}%
which implies that 
\begin{equation*}
\left\Vert w_{\alpha ,\beta ,k}\right\Vert _{L^{2}(\Omega _{\ell
_{o}}^{\prime })}\leqslant \frac{C}{\ell ^{2r}}.
\end{equation*}%
In other way 
\begin{eqnarray*}
\delta ^{\alpha }D^{\beta }\left( u_{\ell }-u_{\infty }\right) 
&\longrightarrow &\frac{1}{\alpha !}D^{\alpha }D^{\beta }\left( u_{\ell
}-u_{\infty }\right) \text{ \ \ in }\mathcal{D}^{\prime }(\Omega _{\ell
_{o}}^{\prime }), \\
\text{\ \ }D^{\beta }\delta ^{\alpha }\left( u_{\ell }-u_{\infty }\right) 
&\longrightarrow &w_{\alpha ,\beta ,k}\text{ \ \ \ \ \ \ in \ }\mathcal{D}%
^{\prime }(\Omega _{\ell _{o}}^{\prime }),
\end{eqnarray*}%
and by uniqueness of the limit, we obtain 
\begin{equation*}
\left\Vert D^{\alpha }D^{\beta }\left( u_{\ell }-u_{\infty }\right)
\right\Vert _{L^{2}(\Omega _{\ell _{o}}^{\prime })}\leqslant \frac{C}{\ell
^{2r}}\text{\ \ for }\left\vert \beta \right\vert \leqslant m.
\end{equation*}%
which gives (\ref{43}), the proof of the theorem is complete.
\end{proof}

\medskip Derivation in the directions $\alpha $ in $N_{1}$ does not get any
trouble to give an estimate on all $\Omega _{\ell _{o}}$, as show it the
following result.

\begin{theorem}
Under assumptions (\ref{3}), (\ref{5}), (\ref{48}), (\ref{27}) and (\ref{45}%
), for any $\ell _{o}>0$ and $r>0$, there exists a constant $C>0$
independent of $\ell $ such that 
\begin{equation*}
\left\Vert D^{\alpha }\left( u_{\ell }-u_{\infty }\right) \right\Vert
_{H^{m}(\Omega _{\ell _{o}})}\leqslant \frac{C}{\ell ^{r}}\text{\ \ for \ }%
\alpha \in N_{1}.
\end{equation*}
\end{theorem}

\begin{proof}
Since (\ref{40}) is verified for $\Omega _{\ell _{o}}^{\prime }=\Omega
_{\ell _{o}}$, $\Omega _{\ell _{1}}^{\prime }=\Omega _{\ell ^{\prime }}$
such that $\ell _{o}<\ell ^{\prime }<\ell _{1}$ for $\alpha \in N_{1}$, then
we can give the same proof as the previous theorem with $\Omega _{\ell
_{o}}^{\prime }=\Omega _{\ell _{o}}.$
\end{proof}

\end{document}